\documentstyle[twoside]{article}

\newtheorem{thm}{\indent{\sc Theorem}}[section]
\newtheorem{defn}[thm]{\indent{\sc Definition}}

\newtheorem{prop}[thm]{\indent{\sc Proposition}} 
\newtheorem{cor}[thm]{\indent{\sc Corollary}}

\newcommand{\ssect}{\subsection}
\newcommand{\sssect}{\subsubsection}

\newcommand{\eqref}[1]{equation~(\ref{#1})}
\newcommand{\eref}[1]{\eqref{#1}}

\newcommand{~}{\nolinebreak[3] }

\newcommand{\th}{\raisebox{0.6ex}{th}}

\newcommand{\seper}{\nopagebreak \begin{center}
	\underline{\hspace{2in}}
	\end{center}}

\newcommand{\bold}[1]{{\em (#1)}\index{#1}}
\newcommand{\bld}[1]{{\em (#1)}}

\newcommand{\dc}[4]{{\samepage$$\left\{ \begin{array}{rrcl}
{\bf Discrete} &\displaystyle #1 &=& \displaystyle  #2 \\*[0.1in]
{\bf Continuous} &\displaystyle #3 &=& \displaystyle  #4 
\end{array}\right\}$$}}

\newcommand{\dclab}[5]{{\samepage \begin{equation}\label{#1}
\left\{ \begin{array}{rrcl}
{\bf Discrete} &\displaystyle #2 &=& \displaystyle  #3 \\*[0.1in]
{\bf Continuous} &\displaystyle #4 &=& \displaystyle  #5
\end{array}\right\}
\end{equation}}}

 \newcommand{\proof}[1]{\proo{#1}$\Box $}
\newcommand{\proo}[1]{{\em Proof: }#1}
\renewcommand{\Box}{ \hspace{1cm}{\rule{1.2ex}{2ex}}}

\newcommand{\rn}[1]{\left\lfloor #1 \right\rceil}

\newcommand{\smallD}{\mbox{{\small D}}}

\newcommand{\D}{{{\rm D}}}

\newcommand{\romcoeff}[2]{\rn{{#1 \atop #2}}}

\newcommand{\action}[2]{\left\langle #2 \right\rangle 
 	\!\!_{_{\scriptstyle #1}}}

\newcommand{\phil}{\nonumber \\*[-1mm]  &&\: }

\newcommand{\seq}[3]{#1_{#2}^{#3}(x)}

\renewcommand{\tilde}[1]{\widetilde{#1}}
  
\newcommand{\fig}[2]{\begin{table}[htbp]
	\caption{#1} 
	\begin{center}\fbox{ \scriptsize \begin{tabular}#2 \end{tabular}}
	\end{center}
\end{table}}

\begin{document}

\newcommand{\thmplus}[2]{\renewcommand{\thethm}{{\bf 
	\thesection{}.\arabic{thm}{\bf #1}}} #2 
	\renewcommand{\thethm}{{\bf \thesection.\arabic{thm}}}}
\renewcommand{\thethm}{{\bf \thesection.\arabic{thm}}}
\newcommand{\explus}[2]{\renewcommand{\theex}{{\bf 
	\thesection{}.\arabic{ex}{\bf #1}}} #2 
	\renewcommand{\theex}{{\bf \thesection.\arabic{ex}}}}
\renewcommand{\theex}{{\bf \thesection.\arabic{ex}}}


\renewcommand{\ssect}[1]{\assectplus{#1}{}{#1}}
\newcommand{\ssectplus}[2]{\assectplus{#1}{#2}{#1}}
\newcommand{\assectplus}[3]{\renewcommand{\thesection}
	{\arabic{section}{\bf #2}}
	\nopagebreak \section[#3]{#1} \nopagebreak \indent}
\renewcommand{\sssect}[1]{\asssectplus{#1}{}{#1}}
\newcommand{\sssectplus}[2]{\asssectplus{#1}{#2}{#1}}
\newcommand{\asssectplus}[3]{\renewcommand{\thesubsection}
	{\thesection.\arabic{subsection}{\bf #2}}
	\nopagebreak \subsection[#3]{#1} \nopagebreak  \indent}
\renewcommand{\numberline}[1]{\makebox[0pt][l]{#1}\hspace{0.8in}}
\renewcommand{\tilde}[1]{\widetilde{#1}}

\begin{abstract}
An extension of the theory of the Iterated Logarithmic Algebra gives the
logarithmic analog of a Sheffer or Appell sequence of polynomials. This leads
to several examples including Stirling's formula and a logarithmic version of
the Euler-MacLaurin summation formula.
\seper 
{\bf L'Alg\`{e}bre des Logarithmes It\'{e}r\'{e}s II: Les Suites de
Sheffer} 
Gr\^ace \`a une g\'en\'eralisation de la th\'eorie de l'alg\`ebre des
logarithmes it\'er\'es, on definit un analogue logarithmique des
suites de polyn\^omes de Sheffer et d'Appell. Quelques exemples
d'applications permettent de d\'{e}duire la formule de Stirling ainsi
qu'un version logarithmique de la formule de sommation de
Euler--MacLaurin.  
 \end{abstract}

\begin{center}
{\it Dedicated to\\
Steve Roman}
\end{center}

\tableofcontents

\ssect{Introduction}
Just as many of the important polynomial sequences of mathematics---generalized
Laguerre, Hermite, Bernoulli---are not quite sequences of binomial type, and
forced us to invent the theory of Sheffer and Appell sequences; we now have in
\cite{ch4} two complete theory of graded sequences of logarithmic series of
binomial type. Both the theories contained in \cite{ch4}---continuous and
discrete---allow us to do work in a field involving the iterated logarithms;
however, the continuous theory also allows us to use $x$ to a real power.
Nevertheless, neither theory can yet help us compute any Euler-MacLaurin-like
summation formulas, since the logarithmic generalization of the Bernoulli
polynomials are not of binomial type. Again, we are forced to augment our
theory through the use of Sheffer and Appell sequences.

Along the way, we will determine a sort of generating function for them, and
for all ordinary graded sequence of formal power series of logarithmic type. We
conclude by presenting several examples and applications. 

The paper \cite{ch4} contains the necessary background for this material.

As in \cite{ch4}, we will study two theories in parallel. All results and
sections regarding the discrete iterated logarithmic algebra will be
denoted ``Discrete'' and the results regarding the more general
continuous iterated logarithmic algebra will be denoted
``Continuous.'' Readers interested in only one of these theories may safely
omit all material pertaining to the other.
Sections and results numbered with a {\bf D} and paragraphs starting with
{\bf Discrete}
are relevant only to the discrete theory whereas sections and results numbered
with a {\bf C} and paragraphs starting with {\bf Continuous} are relevant only
to the continuous theory. Otherwise, the remainder of this article may be
interpreted discretely by 
supposing that all variables $a,b,c,\ldots $ are integers, and that
$\alpha ,\beta $ are vectors of integers, or it may be
interpreted continuously by supposing that all variables
$a,b,c,\ldots $ are real numbers and that $\alpha ,\beta $
are vectors of real numbers.

Digressions and all sections marked with the word ``Appendix'' or
the letter ``{\bf{A}}'' are independent
of all later material, and are included for their own sake.

\ssect{Sheffer Graded Sequences}

We begin by giving one definition of a Sheffer graded sequence of formal power
series of logarithmic type.
\begin{defn}
\bold{Sheffer Graded Sequences} 
Suppose $\seq{s}{a}{\alpha }$ is a graded sequence such that
if there is a Roman graded sequence 
$\seq{p}{a}{\alpha}$ and an Artinian operator $g(\D )$ of
degree zero with $s_{a}^{\alpha }(x)=g(\D )p_{a}^{\alpha
}(x)$. Then we say that $\seq{s}{a}{\alpha }$ is {\em
Sheffer}  for $g(\D )$ with respect to $\seq{p}{a}{\alpha}$.

If $\seq{p}{a}{\alpha}$ is the ($n\th $) associated sequence
for the delta operator $f(\D )$, then we also say that
$\seq{s}{a}{\alpha}$ is {\em Sheffer} for $g(\D )$ with
respect to $f(\D )$ (and $n$).
\end{defn}

Graded Sheffer sequences admit a number of equivalent
characterizations:

\begin{thm}\label{ShChar}
For a graded sequence $\seq{s}{a}{\alpha }$, the following
are equivalent:
\begin{enumerate}
\item $\seq{s}{a}{\alpha}$ is a Sheffer graded sequence.
\item There are Artinian operators $f(\D )$ and $g(\D )$ of
degrees 1 and 0 respectively and an integer $n$ such that 
$$ \action{\alpha }{g(\D )f(\D )^{b;n}s_{a}^{\alpha }(x)} =
\rn{a}!\delta _{ab} $$
for all $a$ and $b$.
\item There is a delta operator $f(\D )$ such that
\begin{equation}
f(\D) s_{a}^{\alpha}(x)=\rn{a}s_{a-1}^{\alpha}(x)
\end{equation}
for all $a$ and $\alpha $.
\item There is a Roman sequence $\seq{p}{a}{\alpha}$ such that
$$ E^{z}s_{a}^{\alpha}(x) =  \sum _{b}
\romcoeff{a}{b} \action{(0)}{E^{z}p_{b}^{(0)}(x)}
s_{a-b}^{\alpha}(x) $$  
for all complex numbers $z$, and all $a$ and $\alpha $.
\end{enumerate}
\end{thm}

\proof{{\bf (4 implies 1)} Let $T$ be the continuous,
linear  operator in ${\cal I}$ 
mapping $p_{a}^{\alpha}(x)$ to $\lambda_{a}^{\alpha}(x)$.
We have: 
\begin{eqnarray*}
TE^{z}s_{a}^{\alpha}(x) &=&\sum_{b} \romcoeff{a}{b}
\action{\alpha} {p_{b}^{\alpha }(x)} Ts_{a-b}^{\alpha}(x)\\* 
&=&\sum_{b} \romcoeff{a}{b}
\action{\alpha} {p_{b}^{\alpha }(x)} Tp_{a-b}^{\alpha}(x)\\
&=&E^{z}p_{a}^{\alpha}(x)\\*
&=&E^{z}Ts_{a}^{\alpha}(x)
\end{eqnarray*} 
 By the characterization of differential operators in \cite{ch4},  $T$ is a
differential operator of degree $0$. Thus, we have
 \begin{equation}\label{T-1}
s_{a}^{\alpha}(x)=T^{-1}p_{a}^{\alpha}(x).
\end{equation}

{\bf (1 implies 3)} We have the following series of
identities: 
\begin{eqnarray*}
f(\D) s_{a}^{\alpha}(x) 
&=& f(\D) T^{-1} p_{a}^{\alpha}(x)\\*
&=& T^{-1} f(\D) p_{a}^{\alpha}(x)\\
&=& \rn{a}T^{-1}p_{a-1}^{\alpha}(x)\\*
&=&\rn{a} s_{a-1}^{\alpha}(x).
\end{eqnarray*}

{\bf (3 implies 4)} Since $E^{z}$ is expressible in term of
$f(\D )^{a;n}$, the result is immediate.

{\bf (2 is equivalent to 1)} See \cite{ch4}.}

\ssect{Expansion Theorem}
As with the harmonic logarithms and all Roman graded sequences, we have
formulas whcih give the coefficients of an arbitrary Artinian operator or
logarithmic series in terms of a Sheffer graded sequence and its operators.
\begin{thm}\bold{Expansion Theorem}
Let the graded sequence $\seq{s}{a}{\alpha}$ be Sheffer for
$g(\D )$ with respect to $f(\D )$ (and $n$). 
 Then for all Artinian operators
$h(\D )$  and  vectors $\alpha \neq  (0)$
we have the following convergent sum.
\begin{equation}\label{GETeq2}
 h(\D ) =\sum_{a} \frac{\action{\alpha}
{h(\D )s_{a}^{\alpha}(x)}}{\rn{a}!} g(\D )^{-1}f(\D )^{a;n}. 
\end{equation}

When $\alpha =(0)$, \eref{GETeq2}holds for all differential operators $h(\D )$.
\end{thm} 

\proof{See the Expansion Theorem of \cite{ch4}.}

\begin{thm}\bold{Taylor's Theorem}\label{GLTF-sh}
Let the graded sequence $\seq{s}{a}{\alpha}$ be Sheffer for
$g(\D )$ with respect to $f(\D )$ (and $n$). 
Then for every formal power series of
logarithmic type $p(x)\in {\cal I}^{+}$ we have the
following convergent sum.
$$ p(x)=\sum_{\alpha\neq (0) }\sum_{a}
{\frac{\action{\alpha}{g(\D)^{-1} f(\D)^{a;n} p(x)}}
{\rn{a}!}} s_{a}^{\alpha}(x).$$ 
\end{thm}

\proof{See the Taylor's Theorem of ~\cite{ch4}.}

\ssect{Generating Functions}\index{Generating Function}

\begin{defn}
\bold{Roman Exponential Series} 
For all $\alpha $,
\begin{description}
\item[Discrete] A vector of integers, define the {\em Roman
exponential series} of $\alpha $ to be the formal series
$$ [e]_{\alpha }^{y}=[exp]_{\alpha }(y)= (E_{\alpha (0)} +
E_{\alpha (1)}) \sum_{n} \frac{y^{n}}{[n]!}. $$
\item[Continuous] A vector of real $\alpha $, define the $n\th $ {\em Roman
exponential series} of $\alpha $ to be the formal series
$$ [e]_{\alpha }^{y;n}=[exp]_{\alpha }(y;n)= (E_{\alpha (0)} +
E_{\alpha (1)}) \sum_{a} \frac{y^{a;n}}{[a]!}. $$
\end{description}
\end{defn}

Thus, the (0\th ) exponential series is a generating
function of $\alpha $ at $x$ is a generating function for
the harmonic logarithms of order $\alpha $.

Note that this is not a formal power series of logarithmic
type in fact is not any sort of Noetherian or Artinian
series; it is  merely a formal series whose coefficients are
homogeneous 
logarithmic series of order $\alpha $. In general, products of
this form are not well defined. For example,
$[e]_{\alpha }^{x}[e]_{\alpha }^{x}$ is not well defined.

In \cite[Corollary 3 to Theorem 2 ]{FOC}, it is shown that
ordinary (nonlogarithmic) 
sequences of binomial type for the operator $f(\D ) $ are
given by the generating function
$$ \sum_{n} \frac{p_{n}(x)}{n!}y^{n} =\exp (xf^{(-1)}(y)).$$
In our present context, we have:

\begin{thm}\index{Generating Function}\label{Roman Generate}
\bld{Generating Function for Roman Graded Sequences} Let 
be  $\seq{p}{a}{\alpha }$ be the ($n\th$) associated sequence of
the delta operator $f(\D )$. Then
\dclab{GenFun}{\sum _{k}
\frac{p_{k}^{\alpha }(x)}{[n]!}y^{n } }{[\exp]_{\alpha}\left(xf^{(-1)}(y)
\right)}{\sum _{a}
\frac{p_{a}^{\alpha }(x)}{[a]!}y^{a}}{[\exp]_{\alpha}^{xf^{(-1)}(y);n}.} 
\end{thm}

\proof{By consideration of
the expansion Theorem \cite{ch4},
$$ p(x,y)=\sum_{a} \frac{\tilde{p}_{a}(x)} {a!} y^{a}
= e^{xf^{(-1;n)}(y);n}. $$  
The result follows now from regularity.}

\begin{thm}\index{Generating Function}
\bld{Generating Function for Sheffer Graded Sequences} Let
$\seq{s}{a}{\alpha }$ be 
Sheffer for $g(\D )$ with respect to $f(\D )$ (and $n$).
Then 
\dc{\sum _{k}
\frac{p_{k}^{\alpha }(x)}{[n]!}y^{n }
}{g(f^{-1}(y))[\exp]_{\alpha}\left(xf^{(-1)}(y) 
\right)}{\sum _{a}
\frac{p_{a}^{\alpha
}(x)}{[a]!}y^{a}}{g(f^{(-1;n)}(y))[\exp]_{\alpha}^{xf^{(-1)}(y);n}.}  
\end{thm}

\proo{By the generalized expansion theorem of \cite{ch4},
\begin{eqnarray*}
\D ^{b} g(\D ) &=& \sum_{a} \frac{\action{\alpha}
{\D ^{b}g(\D )p_{a}^{\alpha}(x)}}{\rn{a}!} f(\D )^{a;n}\\*
f^{(-1;n)}(\D)^{b} g(f^{(-1;n)}(\D )) &=& \sum_{a} \frac{\action{\alpha}
{\D ^{b}g(\D )p_{a}^{\alpha}(x)}}{\rn{a}!} \D ^{a}\\*
&=& \sum_{a} \frac{\action{\alpha}
{\D ^{b} s_{a}^{\alpha}(x)}}{\rn{a}!} \D ^{a}.\Box 
\end{eqnarray*}}

\ssect{Appell Graded Sequences}

The most interesting types of Sheffer graded sequences are
the Roman graded sequences themselves and the Appell graded
sequences. The name {\em Appell} is given to all Sheffer graded
sequences with respect to $\D $ and $0$, or equivalently with
respect to the harmonic logarithms $\lambda_{a}^{\alpha}(x)$. 

They are the  logarithmic generalization of the notion of an Appell
sequence of polynomials, that is, a sequence
$(p_{n}(x))_{n\geq 0}$ of polynomials satisfying the
identity
$$ p_{n}(x+a)=\sum_{k\geq 0}{n\choose k }a^{n-k}p_{k}(x). $$

By \cite{ch4}, the action of the derivative 
operator $\D $ on the space ${\cal I}$ of formal power
series of logarithmic type is naturally  decomposes
as a direct sum of the minimal invariant subspaces ${\cal
I}^{\alpha}$.  For each $\alpha$,
the subspace ${\cal I}^{\alpha}$ is the minimal invariant
subspace of ${\cal I}^{+}$ under the action of the
operators $\D^{a}$ which contains
$\ell^{\alpha}$. It is not yet clear, however, that the
pseudobasis of the spaces ${\cal I}^{\alpha}$ provided by the harmonic
logarithms $\lambda_{a}^{\alpha}(x)$ is determined by
intrinsic algebraic properties. In order to derive the
properties that single out the harmonic logarithms
as the natural pseudobasis for ${\cal I}^{\alpha}$, we are led to
this generalization to formal power series of logarithmic type
of the classical theory of Appell polynomials.

\begin{prop}\label{apTFAE}
Let  $\seq{p}{a}{\alpha}$ be a graded sequence of formal power series of
logarithmic type. Then the following
statements are equivalent:
\begin{enumerate}
\item $\seq{p}{a}{\alpha}$ is an Appell graded sequence.
\item There is an Artinian operator $g(\D )$ of
degree 0  such that
$$ \action{\alpha }{g(\D )\D ^{b} p_{a}^{\alpha }(x)} =
\rn{a}!\delta _{ab} $$
for all $a$ and $b$, and for all $\alpha \neq (0)$.
\item For all $a$ and $\alpha $,
\begin{equation}\label{Dappell} 
\D p_{a}^{\alpha}(x)=\rn{a}p_{a-1}^{\alpha}(x).
\end{equation}
\item For all complex numbers $z$, and all $a$ and $\alpha $,
$$ E^{z} p_{a}^{\alpha}(x) =  \sum _{n\geq 0}
\frac{z^{n}\rn{a}!}{\rn{a-n}!} p_{a-b}^{\alpha}(x). $$  
\end{enumerate}
\end{prop}

\proof{Theorem~\ref{ShChar}.}

We now deduce an explicit expression for an Appell graded sequence as
a linear combination of harmonic logarithms.

\begin{cor}\label{AppellFormula}
Let $\seq{p}{a}{\alpha}$ be an Appell graded sequence, then
for all $a$  and $\alpha $,
$$ p_{a}^{\alpha}(x) =\sum_{b}
\romcoeff{a}{b} \action{(0)}{p_{b}^{(0)}(x)} \lambda
_{a-b}^{\alpha}(x).\Box $$
\end{cor}

\ssect{Examples}
\sssect{Harmonic Graded Sequence}

We begin by giving examples of Appell graded sequences. 
The only Appell graded sequence which is Roman is the graded
sequence of harmonic logarithms. It is
is characterized among 
all Appell graded sequences by the fact that
$$ \action{\alpha } {\lambda _{a}^{\alpha}(x)}=\rn{a}!\delta _{a,0}. $$
This is the characterization of the harmonic graded sequence we had
previously announced.

\sssect{Bernoulli Graded Sequence}
Next, we consider the logarithmic extension of
Bernoulli polynomials.
\begin{defn}\bld{Logarithmic Bernoulli Graded Sequence}
\index{Bernoulli Sequence}\label{BernDef}
Define the {\em Bernoulli operator}\index{Bernoulli Operator} $J$ by
$$ Jp(x)=\int_{x}^{x+1}p(t)dt $$
for $p(x)\in {\cal I}$, that is,
$$J=\frac{e^{\D }-I}{\D }.$$ 
The Appell graded sequence $\seq{B}{a}{\alpha}$ defined as
$$ B_{n}^{\alpha}(x)=J^{-1}\lambda _{n}^{\alpha}(x) ,$$
 is be called the {\em  logarithmic Bernoulli graded sequence}. 
In particular, we obtain the ordinary {\em Bernoulli
polynomials}\index{Bernoulli Polynomials}, $B_{n}(x)=B_{n}^{(0)}(x)$ and the
{\em Bernoulli numbers}\index{Bernoulli Numbers}  
$B_{n}=B_{n}(0)$. 
\end{defn}

The logarithmic Bernoulli graded sequence can be computed by
Proposition~\ref{AppellFormula}:
\begin{equation}\label{Bern*}
 B_{a}^{\alpha}(x)=\sum_{k\geq 0}\romcoeff{n}{k} B_{k} \lambda
_{a-k}^{\alpha}(x). 
\end{equation}

\fig{Logarithmic Bernoulli Graded Sequence
$B_{n}^{\alpha}(x)$}{{rcl}\label{Bern}
	$B_{-2}^{\alpha}(x)$&=&$\lambda _{-2}^{\alpha}(x)+\lambda
_{-3}^{\alpha}(x) +\lambda _{-4}^{\alpha}(x)/2
-\lambda _{-6}^{\alpha}(x)/6 +\lambda_{-8}^{\alpha}(x)/6- \cdots $\\
	$B_{-1}^{\alpha}(x)$&=&$\lambda _{-1}^{\alpha}(x)+
\lambda_{-2}^{\alpha}(x)/2 +\lambda _{-3}^{\alpha}(x)/6
-\lambda_{-5}^{\alpha}(x)/30 +\lambda _{-7} ^{\alpha}(x)/42 -\cdots$\\ 
	$B_{0}^{\alpha}(x)$&=&$\lambda _{0}^{\alpha}(x)
-\lambda_{-1}^{\alpha}(x)/2 -\lambda _{-2}^{\alpha}(x)/12
+ \lambda _{-4}^{\alpha}(x)/120 -\lambda _{-6 }^{\alpha}(x)/252 +\cdots $ \\
	$B_{1}^{\alpha}(x)$&=&$ \lambda _{1}^{\alpha}(x)-
\lambda_{0}^{\alpha}(x)/2 +\lambda _{-1}^{\alpha}(x)/12
-\lambda_{-3}^{\alpha}(x)/360 +\lambda _{-5}^{\alpha}(x)/1260 -\cdots$\\
	$B_{2}^{\alpha}(x)$&=&$\lambda _{2}^{\alpha}(x)-
\lambda_{1}^{\alpha}(x) +\lambda _{0}^{\alpha}(x)/6
+\lambda _{-2}^{\alpha}(x)/360 -\lambda _{-4}^{\alpha}(x)/2520 +\cdots  $}

In particular, 
 the residual series is given by
$$ B_{-1}^{(1)}(x) =
\frac{1}{x} +\frac{1}{2x^{2}} +\frac{1}{6x^{3}}
-\frac{1}{30x^{5}} +\frac{1}{42x^{7}} -\cdots.$$

We derive the Euler-MacLaurin Summation Formula\index{Euler-MacLaurin Summation
Formula} from  the logarithmic Bernoulli graded sequence.  Since $\D 
J=e^{\smallD }-I=\Delta $ where $\Delta $ is the classical forward 
difference operator $\Delta p(x)=p(x+1)-p(x)$,
 the Euler-MacLaurin formula can be written as
$$ I=B_{0}J+B_{1}\Delta +\frac{B_{2}}{2!}\Delta \D
+\frac{B_{3}}{3!}\Delta \D ^{2}+\cdots . $$
Applying it to a discrete formal power series of
logarithmic type $p(x)$ we obtain
\begin{eqnarray*}
\lefteqn{p(x)+p(x+1)+\cdots +p(x+n) = }\\*[-1mm]
 & & 
B_{0}\left[
\int_{x}^{x+n+1}p(t)dt\right]
+B_{1}\left[(p(x+n+1)-p(x)\right]
\phil 
 + \frac{B_{2}}{2!} \left[p'(x+n+1)-p'(x)\right]
+ \frac{B_{3}}{3!} \left[p''(x+n+1)-p''(x)\right] +\cdots .
\end{eqnarray*}
We stress the fact that this formula is an {\it identity} in
the logarithmic algebra, and not just an asymptotic formula. For 
example, for $p(x)=1/x$, we obtain:
\begin{eqnarray}
\lefteqn{ \frac{1}{x}+\frac{1}{x+1}+\cdots +\frac{1}{x+n} = }\nonumber \\*
 & & B_{0}\left[\log (x+n+1)-\log
(x)\right]+B_{1}\left[(x+n+1)^{-1}-x^{-1}\right]
\phil
 +\frac{B_{2}}{2!}\left[- (x+n+1)^{-2}+x^{-2} \right] +
\frac{B_{3}}{3!} \left[2(x+n+1)^{-3}-2x^{-3}\right]+\cdots
.\label{Bern**} 
\end{eqnarray}

For another example, let $p(x)=\log x$. We then obtain a version of Stirling's
formula
\begin{eqnarray}
\lefteqn{\log (x(x+1)\cdots (x+n)) = }\nonumber \\*
 & &  B_{0}((x+n+1)\log
(x+n+1)-x\log x -n-1)
\phil
 +B_{1}(\log (x+n+1)-\log x)
+\frac{B_{2}}{2!}\left[\frac{1}{x+n+1}-
\frac{1}{x}\right] +\cdots . \label{Bern***}
\end{eqnarray}

Finally, let $p(x)=\log \log x$.
\begin{eqnarray*}
\lefteqn{\log (\log x \log (x+1)\cdots \log (x+n) )= }\\*
 & & B_{0}\left[y\log \log y+\sum _{n\geq 1} (n-1)!y (\log
y)^{-n}\right]_{y=x}^{x+n+1}   
+ B_{1}\left[\log \log y \right]_{y=x}^{y+n+1} \\*
&& + \frac{B_{2}}{2}\left[\frac{1}{y\log y}\right]_{y=x}^{x+n+1} 
- \frac{B_{3}}{6}\left[\frac{1+\log y}{y^{2}(\log
y)^{2}}\right]_{y=x}^{x+n+1} +\cdots 
\end{eqnarray*}
We note that $J=\Delta \D ^{-1}$, and thus for $\alpha \neq (0)$,
\begin{eqnarray*}
\Delta B_{a}^{\alpha}(x)&=&\Delta J^{-1}\lambda _{a}^{\alpha}(x)\\*
&=&\D \lambda_{a}^{\alpha}(x)\\*
&=&\rn{a}\lambda _{a-1}^{\alpha}(x).
\end{eqnarray*}
Summing, we obtain
$$ \lambda _{a}^{\alpha}(x)+\lambda _{a}^{\alpha}(x+1)+\cdots
+\lambda _{a}^{\alpha}(x+k)=\rn{a} ^{-1}\left[B_{a+1}^{\alpha}(x+k+1)
-B_{a+1}^{\alpha}(x)\right]. $$            
For example, 
\begin{eqnarray*}
\log (x(x+1)\cdots (x+k)) 
&=& B_{1}^{(1)}(x+k+1)-B_{1}^{(1)}(x)\\*
\log (\log x\log (x+1)\cdots \log (x+k)) 
&=& B_{1}^{(0,1)}(x+k+1)-B_{1}^{(0,1)}(x)\\
(\log x)^{2}+\cdots +(\log (x+n+1))^{2} &=&
B_{1}^{(2)}(x+k+1)-B_{1}^{(2)}(x)\\*
\end{eqnarray*}
Clearly, any similar sum of $\ell ^{(a),\alpha }$ can be
evaluated in terms of the 
logarithmic Bernoulli graded sequence.

\sssect{Hermite Graded Sequence}

Our next example of an Appell graded sequence is the logarithmic
Hermite graded sequence.

\begin{defn}
\bld{Logarithmic Hermite Graded Sequence}\index{Hermite Sequence}
Let the {\em Weirstrass operator}\index{Weirstrass Operator} be $W=e^{\smallD
^{2}/2}$. The {\em logarithmic Hermite graded sequence}
$\seq{H}{a}{\alpha}$ is defined as
$$ H_{a}^{\alpha}(x)= W^{-1}\lambda _{a}^{\alpha}(x). $$
\end{defn}

For $n$ a nonnegative integer,
$H_{n}^{(0)}(x)=H_{n}(x)$ is the usual 
Hermite polynomial\index{Hermite Polynomial}. From them we define the {\em
Hermite numbers}\index{Hermite Numbers}The logarithmic Hermite graded
sequence 
can be computed via Proposition~\ref{AppellFormula}:
\begin{equation}\label{correction}
H_{a}^{\alpha}(x)=\sum_{k\geq 0}\romcoeff{a}{k}
H_{k} \lambda _{a-k}^{\alpha}(x).
\end{equation}
\fig{Logarithmic Hermite Graded Sequence
$H_{n}^{\alpha}(x)$}{{rcl}\label{Herm} 
	$H_{-2}^{\alpha}(x)$&=&
$\lambda _{-2}^{\alpha}(x)
-6\lambda_{-4}^{\alpha}(x)
+60\lambda _{-6}^{\alpha}(x)
-840\lambda_{-8}^{\alpha}(x)- \cdots $\\
	$H_{-1}^{\alpha}(x)$&=&
$\lambda _{-1}^{\alpha}(x)
-2\lambda _{-3}^{\alpha}(x)
+12\lambda_{-5}^{\alpha}(x)
-120\lambda _{-7} ^{\alpha}(x)-\cdots $\\
	$H_{0}^{\alpha}(x)$&=&
$\lambda _{0}^{\alpha}(x) 
+\lambda _{-2}^{\alpha}(x)
-3\lambda _{-4}^{\alpha}(x)
+20\lambda _{-6 }^{\alpha}(x)+\cdots $\\
	$H_{1}^{\alpha}(x)$&=&
$ \lambda _{1}^{\alpha}(x)
-\lambda _{-1}^{\alpha}(x)
+\lambda_{-3}^{\alpha}(x)
-4\lambda _{-5}^{\alpha}(x)-\cdots$\\
	$H_{2}^{\alpha}(x)$&=&
$\lambda _{2}^{\alpha}(x)
-2\lambda _{0}^{\alpha}(x)
-\lambda _{-2}^{\alpha}(x)
+2\lambda _{-4}^{\alpha}(x)+\cdots  $}
(See Table~\ref{Herm}.) Recall that
$$ E_{0\alpha }H_{n}^{\alpha}(x) =H_{n}(x),$$
in other words, the ``positive'' terms of \eref{correction}
equal the ordinary Hermite polynomials, except that $\lambda_{n}^{\alpha}(x)$
replaces $x^{n}$.  Thus, we have  a 
logarithmic generalization of the Hermite polynomials.  All
classical properties of Hermite 
polynomials may be extended to the logarithmic Hermite
graded sequence. For example, we have the  explicit expression
\begin{eqnarray}
H_{a}^{\alpha}(x) &=& e^{-\smallD ^{2}/2}\lambda _{a}^{\alpha}(x)\nonumber \\*
&=&\sum_{k\geq
0}\left(-\frac{1}{2}\right)^{k}\frac{\rn{a}!}{k!\rn{a-2k}!}
\lambda _{a-2k}^{\alpha}(x). \label{H1}
\end{eqnarray}
so that for $n$ a negative integer,
\begin{equation}\label{H2}
H_{n}^{(1)}(x)=\sum_{k\geq
0}\left(-\frac{1}{2}\right)^{k}\frac{(2k-n-1)!}{k!(-n-1)!}x^{n-2k}
\end{equation}
or
$$ H_{n}^{\alpha }(x)=\sum_{k\geq 0}
\left(-\frac{1}{2}\right)^{k} \frac{(2k-n-1)!}{k!(-n-1)!}
\lambda_{n-2k}^{\alpha}(x)
 $$
The right side of \eref{H2} is the classical asymptotic expansion of the 
Hermite series $H_{n}^{(1)}(x)$; in the present context, it
is  a convergent series, and one term of the Hermite graded
sequence.  We trivially have
$$ \D H_{a}^{\alpha}(x)=\rn{a}H_{a-1}^{\alpha}(x), $$
and
$$ E^{z}H_{a}^{\alpha}(x)=\sum_{k\geq 0}\romcoeff{a}{k}
z^{k}H_{a-k}^{\alpha}(x).$$
Finally, every logarithmic series $p(x)$ has a unique convergent
expansion in terms of the logarithmic Hermite graded sequence:
$$ p(x) = \sum_{a,\alpha }\frac{\action{\alpha } {W\D ^{a}p(x)}}{\rn{a}!}
H_{a}^{\alpha}(x).$$ 

\sssect{Laguerre Graded Sequence}

We conclude by introducing a Sheffer sequence which is
neither Appell nor Roman. The generalized Laguerre graded sequence is Sheffer
for the ordinary Laguerre graded sequence \cite{ch4}.

\begin{defn} 
\bld{Laguerre Graded Sequence of Grade $b$}
\index{Laguerre Sequence}
Given a real number $b$, define the {\em
Laguerre sequence of grade $b$} $L_{a}^{\alpha ;b}(x)$ to
be the Sheffer graded sequence for $(1-\D )^{b+1;0}$ with
respect to the Laguerre graded sequence.
\end{defn}

By the transfer formula \cite{ch4},
\begin{eqnarray*}
 L_{a}^{\alpha ;b}(x) &=& (1-\D )^{b+1;0} K' \left(\frac{\D
}{K} \right)^{a+1;0} \lambda_{a+1}^{\alpha}(x)\\*
&=& (-1)^{a;0} (1-\D )^{a+b} \lambda_{a}^{\alpha}(x)\\*
&=& \sum _{k\geq 0} {a+b \choose k}
\frac{\rn{a}!}{\rn{a-k}!} (-1)^{a-k;0}
\lambda_{a-k}^{\alpha}(x) 
\end{eqnarray*}

The generating function for the Laguerre graded sequence of
grade $b$ is 
$$ (1-y)^{-b-1;0} \exp _{\alpha } (xy/(y-1)) =\sum _{a}
\frac{L_{a}^{\alpha ;b}}{\rn{a}!}y^{a}.$$

We hope the preceding examples display the utility of the
theory of formal power series of logarithmic type.

Please refer to \cite{ch4} for a more complete list of references on this
topic. 


\begin{thebibliography}{99}
\bibitem{ch4}
{\sc D. Loeb}, {\em The Iterated Logarithmic Algebra}, To appear.
\bibitem{33}
{\sc J. B. Miller},
{\em The Standard Summation Operator, the Euler-MacLaurin
Sum Formula, and the Laplace Transformation}, Journal of the
Australian Mathematical Society,
{\bf 39} (1985), 376--390.
\bibitem{37}
{\sc T. R. Prabhakar, and Reva},
{\em An Appell Cross-sequence Suggested by the Bernoulli and
Euler Polynomials of General Order}, Indian Journal of Pure
and Applied Mathematics,
{\bf 10} (1979), 1216--1227.
\bibitem{41}
{\sc N. Ray},
{\em Extensions of Umbral Calculus, Penumbral Coalgebras
and Generalized Bernoulli Numbers}, Advances in Mathematics, 
{\bf 61} (1986), 49--100.
\bibitem{FOC}
{\sc G.-C. Rota}, ``Finite Operator Calculus,'' Academic
Press, 1975.
\end{thebibliography}
\end{document}